\documentclass[12pt,twoside]{amsart}
\usepackage{amssymb}

\nonstopmode

\numberwithin{equation}{section}
\font\tengothic=eufm10 scaled\magstep 1
\font\sevengothic=eufm7 scaled\magstep 1
\newfam\gothicfam
\textfont\gothicfam=\tengothic
\scriptfont\gothicfam=\sevengothic

\newcommand{\fc}{\mathfrak c}
\newcommand{\Z}{\mathbb{Z}}

\DeclareMathOperator{\Tor}{Tor}

\newcommand{\ffi}{\varphi}
\newcommand{\al}{\alpha}
\newcommand{\be}{\beta}

\DeclareMathOperator{\pnt}{\raise 0.5mm \hbox{\large\bf.}}

\newtheorem{theorem}{Theorem}[section]

\newtheorem{proposition}[theorem]{Proposition}
\newtheorem{corollary}[theorem]{Corollary}
\newtheorem{conjecture}[theorem]{Conjecture}
\newtheorem{claim}[theorem]{Claim}

\theoremstyle{definition}
\newtheorem{remark}[theorem]{Remark}
\newtheorem{example}[theorem]{Example}

\begin{document}

\title{The multiplicitiy conjecture in low codimensions}
\author{Juan Migliore}
\address{Department of Mathematics, University of Notre Dame,
Notre Dame, IN 46556, USA}
\email{Juan.C.Migliore.1@nd.edu}
\author{Uwe Nagel}
\address
{Department of Mathematics, University of Kentucky,
715 Patterson Office Tower, Lexington, KY 40506-0027, USA}
\email{uwenagel@ms.uky.edu}
\author{Tim R\"omer}
\address{FB Mathematik/Informatik,
Universit\"at Osnabr\"uck, 49069 Osnabr\"uck, Germany}
\email{troemer@mathematik.uni-osnabrueck.de}

\thanks{Part of the work for this  paper was done while the first author 
was sponsored by the National Security Agency  under Grant Number
MDA904-03-1-0071 and the second author was supported by
a Special Faculty Research Fellowship from the University of Kentucky.}

\begin{abstract}
We establish the multiplicity conjecture of Herzog, Huneke, and Srinivasan 
about the multiplicity of graded Cohen-Macaulay algebras over a field, for
codimension two algebras and for Gorenstein algebras of codimension
three. In fact, we prove stronger bounds than the conjectured ones
allowing us to characterize the extremal cases. This may be seen as a
converse to the multiplicity formula of Huneke and Miller that inspired
the conjectural bounds.
\end{abstract}

\maketitle
\tableofcontents

\section{Introduction}
Let $R = K[x_1,\dots,x_n]$ be
the polynomial ring over a field $K$
with its standard grading where $\deg x_i=1$ for $i=1,\dots,n$.
Let $R/I$ be a standard graded $K$-algebra,
where $I$ is a graded ideal of codimension $c$.
We denote by $e(R/I)$ the multiplicity of $R/I$.
When $I$ is a saturated ideal defining a closed subscheme
$V \subset {\mathbb P}^{n-1}$,
$e(R/I)$ is just the degree $\deg V$ of $V$.
Consider the minimal graded free resolution of $R/I$
$$
0 \rightarrow
\bigoplus_{j \in \Z} R(-j)^{\beta_{p,j}^R(R/I)}
\rightarrow \dots \rightarrow
\bigoplus_{j \in \Z} R(-j)^{\beta_{1,j}^R(R/I)}
\rightarrow R
\rightarrow 0
$$
where
we denote by
$\beta_{i,j}^R(R/I)=\Tor_i^R(R/I,K)_j$ the graded Betti numbers of $R/I$
and $p$ is the projective dimension of $R/I$.
Let $c$ denote the codimension of $R/I$.
Then $c \leq p$ and equality holds if and only if
$R$ is Cohen-Macaulay.
We define
$m_i = \min \{j \in \Z :  \beta_{i,j}^R(R/I)\neq 0 \}$ and
$M_i = \max \{j \in \Z :  \beta_{i,j}^R(R/I)\neq 0 \}$.
When there is any danger of ambiguity, we write
$m_i(I)$ and $M_i(I)$.
The algebra
$R/I$ has a {\em pure resolution} if
$m_i = M_i$ for all
$i$.
In this case we write $d_i$ for the
unique $j$ such that $\beta_{i,j}^R(R/I)\neq 0$.
It was shown by Huneke and Miller \cite{HM}
that if $R/I$ is Cohen-Macaulay with a pure resolution
then
\[
e(R/I) = \left ( \prod_{i=1}^p d_i \right )/ p! .
\]
Extending this, Herzog, Huneke and Srinivasan made the following
multiplicity conjecture:

\begin{conjecture} \label{conj}
If $R/I$ is Cohen-Macaulay then
$$
\left ( \prod_{i=1}^p m_i \right )/ p! \leq e(R/I) \leq \left (
\prod_{i=1}^p M_i \right )/p!.
$$
\end{conjecture}
Conjecture \ref{conj} has been studied extensively, and partial results have
been obtained.  In \cite{HS}, Herzog and Srinivasan proved it
in the following cases:
$R/I$ is a complete intersection;
$I$ is a perfect ideal with quasi-pure resolution (i.e. $m_i(R/I)\geq M_{i-1}(R/I)$);
$I$ is a perfect ideal of codimension 2;
$I$ is a codimension 3 Gorenstein ideal with five minimal generators;
$I$ is a Gorenstein monomial ideal of codimension 3 with at
least one generator of smallest possible degree (relative to the number of generators);
$I$ is a perfect stable ideal;
$I$ is a perfect squarefree strongly stable ideal.
Furthermore,
Herzog and Srinivasan proved that the {\em upper bound} of Conjecture \ref{conj} holds for
all codimension 3 Gorenstein ideals.
In addition, Guardo and Van Tuyl
\cite{GV} proved that the conjecture holds for powers of complete
intersections, and Gold, Schenck and Srinivasan
\cite{GSS} proved it in certain cases where $I$ is linked to something
``simple."  In addition, Srinivasan \cite{srinivasan}
proved a stronger bound for Gorenstein algebras with quasi-pure resolutions.
(cf.\ Remark \ref{rem-hema}.)

The non Cohen-Macaulay case has also been studied.
Here it is necessary to replace the projective dimension by the codimension in
Conjecture \ref{conj}.  It was observed in \cite{HS} that the lower bound is
false.  However, Herzog and Srinivasan proved the upper bound in these cases:
$I$ is a stable ideal;
$I$ is a squarefree strongly stable ideal;
$I$ is an ideal with $d$-linear resolution.
In addition, Gold \cite{gold} proved it for codimension two lattice ideals; this was
generalized by R\"omer \cite{roemer} for all codimension two ideals.
It was also proved by Gasharov, Hibi and Peeva \cite{GHP}
for {\bf a}-stable ideals
and more generally by R\"omer \cite{roemer} for componentwise linear ideals.

In this paper we begin with a new, stronger version of Conjecture \ref{conj} in
the codimension two case:

\begin{theorem} \label{thm-intro-2}
Let $R/I$ be a graded Cohen-Macaulay algebra of codimension two.
Then the following lower and upper bounds hold:
\begin{enumerate}
\item
$
e(R/I)
\geq
\frac{1}{2} m_1 m_2 + \frac{1}{2} (M_2-M_1) (M_2-m_2 + M_1-m_1)
$;
\item
$
e(R/I) \leq \frac{1}{2} M_1 M_2 -  \frac{1}{2} (m_2-m_1)
( M_2-m_2 + M_1-m_1 ).
$
\end{enumerate}
\end{theorem}

As an immediate consequence of this result, we obtain the following
characterization for the sharpness of Conjecture \ref{conj} in the codimension two
  Cohen-Macaulay case.  This can be viewed as a converse to the Huneke-Miller
result \cite{HM} mentioned above.

\begin{corollary} \label{cor-intro-2}
Let $R/I$ be a graded Cohen-Macaulay algebra of codimension two. Then the following conditions
are equivalent:
\begin{enumerate}
\item $e(R/I) = \frac{1}{2} m_1 m_2 $;
\item $e(R/I) = \frac{1}{2} M_1 M_2 $;
\item $R/I$ has a pure minimal graded free resolution.
\end{enumerate}
\end{corollary}

In Section \ref{sec-2} we give two proofs of these results.
The first is based on some formulas in \cite{HS}.
The second one is more self-contained and uses the methods
that allow us obtain the results discussed below.
We also discuss the non Cohen-Macaulay case.
In particular, we show that even a natural
weakening of the lower bound in Conjecture \ref{conj}
could only  be true for reduced ideals (cf.\ Remark \ref{rem-naCM}).
In Section \ref{sec-3} we prove a stronger version of
Conjecture \ref{conj} for Gorenstein ideals of codimension three:

\begin{theorem} \label{ht 3 gor bd}
Let $R/I$ be a graded Gorenstein algebra of codimension three.
Then the following lower and upper bounds hold:
\begin{enumerate}\label{eq-gor}
\item
$e(R/I)\geq \frac{1}{6} m_1 m_2 m_3 + \frac{1}{6}(M_3-M_2)2(M_2-m_2+M_1-m_1)$;
\item
$e(R/I) \leq \frac{1}{6} M_1 M_2 M_3 - \frac{1}{12} M_3(M_2-m_2+M_1-m_1)$.
\end{enumerate}
\end{theorem}

As an immediate consequence we get
that the lower bound of Conjecture \ref{conj} holds
for codimension three Gorenstein ideals.
This was the last open case of the conjecture
in low codimensions where structure theorems of
minimal graded free resolutions are available.
As in the case of codimension two perfect ideals we can characterize
when Conjecture \ref{conj} is sharp providing again a converse to the Huneke-Miller
formula in  \cite{HM}.

\begin{corollary} \label{cor-sharp-gor}
Let $R/I$ be a graded Gorenstein algebra of codimension three.
Then the following
conditions are equivalent.
\begin{enumerate}
\item $e(R/I) = \frac{1}{6} m_1 m_2 m_3$;
\item $e(R/I) = \frac{1}{6} M_1 M_2 M_3 $;
\item $R/I$ has a pure minimal free resolution.
\end{enumerate}
\end{corollary}

Our method of proof consists in exhibiting
a specific example for each possible set of Betti numbers and a reduction procedure
that allows us to proceed by induction.  While the first idea seems difficult to
extend we expect that the reductions via basic double links will be useful in other cases, too.

We conclude this note with an explicit formula for the multiplicity of a Gorenstein ideal
in terms of the degrees of the entries of its Buchsbaum-Eisenbud matrix (Proposition
\ref{prop-def-formula}). It is a by-product of the proof of Theorem \ref{ht 3 gor bd}.

\section{Codimension two Cohen-Macaulay algebras}

\label{sec-2}
Let $K$ be a field,
$R=K[x_1,\dots,x_n]$ be the polynomial ring and
$I \subset R$ a graded ideal of height two
such that $R/I$ is a Cohen-Macaulay ring.
It follows from the Hilbert-Burch theorem (e.g. see \cite{BRHE} for details)
that $I$ has a minimal graded free resolution of the form
\begin{equation} \label{eq-2-res}
0 \to
\bigoplus_{i=1}^{m-1}R(-f_i)
\stackrel{\ffi}{\longrightarrow}
\bigoplus_{i=1}^{m}R(-e_i)
\to
J
\to
0
\end{equation}
Let $u_i=f_i-e_i$ and $v_i=f_i-e_{i+1}$.
The following was observed, for example, in \cite{HTV}:
\begin{enumerate}
\item
$u_i \geq v_i \geq 0$ for $i=1,\dots,m-1$;
\item
$u_{i+1} \geq v_i$ for $i=1,\dots,m-2$;
\item
$e_1=v_1+\dots+v_{m-1}$
and
$e_m=u_1+\dots+u_{m-1}$;
\item
$f_1=v_1+\dots+v_{m-1}+u_1$
and
$f_{m-1}=u_1+\dots+u_{m-1}+v_{m-1}$;
\item
$
e(R/I)=
\sum_{i=1}^{m-1} u_i(v_i+\dots + v_{m-1})
=
\sum_{i=1}^{m-1} v_i(u_1+\dots + u_{i}).
$
\end{enumerate}

Herzog and Srinivasan proved in
\cite{HS}
the formulas (see proof of Theorem 2.1):
\begin{enumerate}
\item
$\sum_{i=2}^{m-1}(v_{i-1} + v_i)(v_i+\dots + v_{m-1})
=
(v_1+\dots+v_{m-1})(v_2+\dots+v_{m-1})
$;
\item
$\sum_{i=1}^{m-2}(u_{i} + u_{i+1})(u_1+\dots + u_{i})
=
(u_1+\dots+u_{m-1})
(u_1+\dots+u_{m-2})
$.
\end{enumerate}

Note that  $e_1=m_1$, $e_{m}=M_1$,
$f_1=m_2$ and $f_{m-1}=M_2$.
Following the proof of  Theorem 2.1 in \cite{HS},
we can show Theorem \ref{thm-intro-2} of the introduction.

\begin{proof}[Proof of Theorem \ref{thm-intro-2}]
(a)
Using
\begin{eqnarray*}
2u_1
=
(u_1-v_1) + u_1+v_1
=
e_2-e_1 + u_1+ v_1
\\
2u_i
=
(2u_i-v_{i-1}-v_i) + v_{i-1}+v_i
=
(e_{i+1}-e_i + f_i - f_{i-1}) + v_{i-1}+v_i
\end{eqnarray*}
and the first formula of Herzog and Srinivasan above
we compute
\begin{eqnarray*}
2e(R/i)
& = &
\sum_{i=1}^{m-1} 2u_i(v_i+\dots + v_{m-1})\\
& = &
e_1f_1
+
(u_1-v_1)(v_1+\dots + v_{m-1})
+
\sum_{i=2}^{m-1} (2u_i-v_i-v_{i-1})(v_i+\dots + v_{m-1})
\\
& = &
m_1m_2
+
(e_2-e_1)(v_1+\dots + v_{m-1})
+
\sum_{i=2}^{m-1} (e_{i+1}-e_i + f_i - f_{i-1})(v_i+\dots + v_{m-1})
\\
& = &
m_1m_2
+
\sum_{i=1}^{m-1}(e_{i+1}-e_i)(v_i+\dots + v_{m-1})
+
\sum_{i=2}^{m-1} (f_i - f_{i-1})(v_i+\dots + v_{m-1})
\\
& \geq &
m_1m_2
+
v_{m-1}
(
\sum_{i=1}^{m-1}(e_{i+1}-e_i)
+
\sum_{i=2}^{m-1} (f_i - f_{i-1})
)
\\
& = &
m_1m_2
+
(M_2-M_1)
(
M_1-m_1
+
M_2-m_2
).
\end{eqnarray*}
Dividing by 2, the desired formula follows.

(b)
Similarly, using the second formula of
Herzog and Srinivasan
we compute
\begin{eqnarray*}
2e(R/I)
& = &
\sum_{i=1}^{m-1} 2v_i(u_1+\dots + u_{i})
\\
& = &
e_m f_{m-1}
+
(v_{m-1}-u_{m-1})(u_1+\dots + u_{m-1})
+
\sum_{i=1}^{m-2} (2v_i-u_i-u_{i+1})(u_1+\dots + u_{i})
\\
& = &
M_1M_2
-
(e_{m}-e_{m-1})(u_1+\dots + u_{m-1})
-
\sum_{i=1}^{m-2} (e_{i+1}-e_i + f_{i+1} - f_{i})(u_1+\dots + u_{i})
\\
& = &
M_1M_2
-
\sum_{i=1}^{m-1}
(e_{i+1}-e_i)(u_1+\dots + u_{i})
-
\sum_{i=1}^{m-2} (f_{i+1} - f_{i})(u_1+\dots + u_{i})
\\
& \leq &
M_1M_2
-
u_1
(
\sum_{i=1}^{m-1}
(e_{i+1}-e_i)
+
\sum_{i=1}^{m-2} (f_{i+1} - f_{i})
)
\\
& = &
M_1M_2
-
(m_2-m_1)
(
M_1-m_1
+
M_2-m_2
)
\end{eqnarray*}
\end{proof}

Now we present  alternative, more self-contained proofs for the bounds in Theorem
\ref{thm-intro-2}  that use some methods of liaison theory. Its purpose is twofold.
They illustrate the principles we will use in the following section to discuss Gorenstein
ideals of
codimension three and they allow us to provide some of the relations we will use there.

The map $\ffi$ in the minimal free resolution (\ref{eq-2-res})
is represented by the Hilbert-Burch matrix of $I$.
By reordering  we can arrange that the degree matrix is
\begin{equation} \label{deg mat}
A = \begin{bmatrix}
a_{1,1} & a_{1,2} & \dots & a_{1,t+1} \\
\vdots & \vdots &  & \vdots \\
a_{t,1} & a_{t,2} & \dots & a_{t,t+1}
\end{bmatrix}
\end{equation}
where the entries are increasing from bottom to top and from left to right, so
$a_{t,1}$ is the smallest and $a_{1,t+1}$ is the largest.
Notice that in order to be a degree matrix of a height two
Cohen-Macaulay ideal,
the main diagonal has to be strictly positive
(\cite{GM3}, page 3142 and see \cite{sauer}, page 84).

\begin{remark} \label{rem-deg-matrix}
Let us rewrite the degree matrix $A$ as follows
$$
A = \begin{bmatrix}
a_1 & b_1 & & *\\
& \ddots & \ddots &  \\
* & & a_t & b_t \\
\end{bmatrix}
$$
Note that $A$ is completely determined by
$a_1,\ldots,a_{t},b_1,\ldots,b_{t}$.
By our ordering of degrees we have in particular
$$
b_t \geq a_t \quad \mbox{and} \quad b_{t-1} \geq a_t \; \mbox{provided} \; t
\geq 2.
$$
The minimal generators of $I$ have degrees
$a_1 + \ldots + a_j + b_{j+1} + \ldots + b_t$, $j = 0,\ldots,t$
and the syzygies of $I$ have degrees
$a_1 + \ldots + a_j + b_{j} + \ldots + b_t$, $j = 1,\ldots,t$.
Thus, we obtain
\begin{equation} \label{formulas}
\begin{array}{rcl}
m_1 & = & a_1 + \ldots + a_t \\
M_1 & = & b_{1} + \ldots + b_t \\
m_2 & = & a_1 + \ldots + a_t + b_t = m_1 + b_t \\
M_2 & = & a_1 + b_{1} + \ldots + b_t.
\end{array}
\end{equation}
\end{remark}

Now we will re-prove the lower bound  for the multiplicity in
Theorem \ref{thm-intro-2}:

\begin{proof}[Proof of Theorem \ref{thm-intro-2}(a)]
We induct on $t \geq 1$. It $t = 1$ then $I$ is a complete intersection. Its
resolution is
\begin{equation} \label{koszul}
0 \rightarrow R(-m_1-M_1) \rightarrow
\begin{array}{c}
R(-M_1) \\
\oplus \\ R(-m_1)
\end{array}
\rightarrow I \rightarrow 0.
\end{equation}
Then we get
\[
\begin{array}{rcl}
\frac{1}{2} m_1 m_2 + \frac{1}{2} (M_2-M_1) (M_2-m_2 + M_1-m_1) &
= & \frac{1}{2} m_1 (m_1+M_1) + \frac{1}{2} m_1(M_1-m_1) \\[1ex]
& = & m_1 M_1 = e(R/I).
\end{array}
\]

Now assume $t \geq 1$ and let $I'$ be an ideal whose degree matrix of the
Hilbert-Burch matrix is
$$
A' = \begin{bmatrix}
a_1 & b_1 & & & *\\
& \ddots & \ddots &  \\
 & & a_t & b_t \\
*   & & & a_{t+1} & b_{t+1}
\end{bmatrix}.
$$
Note
that the multiplicity of $R/I'$ is completely determined by the degree matrix,
so it suffices to consider an example of an ideal for any degree matrix.  Basic double
linkage is then used in order to apply the induction hypothesis.  The idea will
be to show that we can reduce to an ideal with degree matrix $A$ (see Remark
\ref{rem-deg-matrix}), i.e.\ we remove the last row and the last
column.

It is easy to see that the following  monomial ideal has $A'$ as its
degree matrix
$$
I' = (y^{b_1 + \ldots + b_{t+1}}, x^{a_1} y^{b_2 + \ldots + b_{t+1}},
\ldots, x^{a_1 + \ldots + a_{t+1}}).
$$
Write it as
\begin{equation} \label{eq-J'-J}
I' = x^{a_1 + \ldots + a_{t+1}}R + y^{b_{t+1}} I.
\end{equation}
Then the monomial ideal $I$ has $A$ as its degree matrix where $A$ is obtained
by deleting the last row and column of $B$. Thus, we may apply induction to
$I$.
Let $m_1,m_2,M_1,M_2$ be the corresponding invariants for $I$, and let
$m_1',m_2',M_1',M_2'$ be those of $I'$. Moreover, in order to simplify notation
we set
$$
a := a_{t+1}, \; b:= b_{t+1}, \; \mbox{and} \; c := b_t.
$$
Using Remark \ref{rem-deg-matrix}, we see that
\begin{equation} \label{comp-deg}
\begin{array}{rcl}
m_1' & = & m_1 + a \\
M_1' & = & M_1 + b \\
m_2' & = & m_2 + a + b -c \\
M_2' & = & M_2 + b.
\end{array}
\end{equation}
Moreover, we have $e(R/I') = e(R/I) + m_1' b$.
Using the formulas above we get
by induction
\begin{eqnarray*}
\lefteqn{ m_1'm_2' +  (M_2'-M_1') (M_2'-m_2' +  M_1'-m_1')} \\[1ex]
   & = &  (m_1+a)(m_2+a+b-c) \\
   & & + [(M_2+b) - (M_1+b)] [(M_2+b) - (m_2+a+b-c) + (M_1 + b) - (m_1+a)] \\[1ex]
& = &  m_1 m_2  +  (m_1 + a) (a+b-c) +  m_2 a \\
& & +  (M_2 - M_1) [M_2 - m_2+ M_1  - m_1 + b + c - 2a] \\[1ex]
& \leq & 2 e(R/I) +  (m_1 + a) (a+b-c) +  m_2 a +
   (M_2 - M_1) (b+c-2a) \\[1ex]
& = & 2 e(R/I') - 2 (m_1 + a) b +  (m_1 + a) (a+b-c) + (m_1 + c)  a + \\
& & (M_2 - M_1) (b+c-2a)
\\[1ex]
 & = & 2 e(R/I') +  [m_1 + M_1 - M_2] (2a-b-c) +  a(a-b).
\end{eqnarray*}
Note that we used the relation $m_2 = m_1 + c$.
But $b \geq a$ and $c \geq a$ and $M_1 + m_1 \geq M_2$ by \ref{rem-deg-matrix}, so
\[
\begin{array}{rcll}
\frac{1}{2} m_1'm_2' + \frac{1}{2} (M_2'-M_1') (M_2'-m_2' +
M_1'-m_1') & \leq & e(R/I')
\end{array}
\]
as desired.  (Note the strict inequality unless $a = b = c$ or $t = 2$ and
$a=b$.)
\end{proof}

In a similar way we can re-prove the upper bound in Theorem
\ref{thm-intro-2}.

\begin{proof}[Proof of Theorem \ref{thm-intro-2}(b)]
We use again induction on $t \geq 1$.  First suppose that $I$
is a complete
intersection.
Then we have the resolution (\ref{koszul}),
and we obtain
\[
\begin{array}{rcl}
\frac{1}{2} M_1 M_2 -
\frac{1}{2}
(m_2-m_1)
(
M_2-m_2
+
M_1-m_1
) & = &
\frac{1}{2} M_1 (M_1 + m_1) -
\frac{1}{2} M_1 (M_1-m_1) \\[1ex]
& = & m_1 M_1 = e(R/I).
\end{array}
\]

Now for the general case, we again use induction with the set-up of
the previous
proof.
We have
\begin{eqnarray*}
\lefteqn{  M_1'M_2' - (m_2'-m_1')
(
M_2'-m_2'
+
M_1'-m_1'
)}
\\[1ex]
 & = &  (M_1+b)(M_2+b)
  \\
& & -  [(m_2 + a + b - c) - (m_1+a)]  [M_2 -
m_2+ M_1  - m_1 + b + c - 2a]
\\[1ex] & = & M_1M_2 +  b (M_1 + M_2)  +  b2 \\
& & -  [m_2 + b - c - m_1]  [M_2 - m_2+ M_1  - m_1 + b + c - 2a]
\\[1ex]
 & \geq & 2 e(R/I) +   b (M_1 + M_2 + b)  -  b
(b+c -2a)  -  (b-c) [M_2 - m_2+ M_1  - m_1]
\\[1ex]
& = & 2 e(R/I') - 2 (m_1 + a) b +  b (M_1 + M_2 + b)  -  b (b+c -2a)  \\
& & - (b-c) [M_2 - m_2+ M_1  - m_1] \\[1ex]
& = & 2 e(R/I') +   b (M_1 + M_2 - 2m_1 -c)   - (b-c) [M_2 - m_2+ M_1 - m_1] \\[1ex]
& = & 2 e(R/I') +   c (M_1 + M_2 - 2m_1 -c)\\[1ex]
& \geq & 2 e(R/I').
\end{eqnarray*}
(Note the strict inequality unless $I$ has a pure resolution.)
\end{proof}

In a similar way we can prove another upper bound
that extends the upper bound of Herzog and Srinivasan.
The following proposition has a
hypothesis that is a bit technical, but it has a more satisfying
conclusion
than the bound of Theorem \ref{thm-intro-2}(b) in one case.

\begin{proposition} \label{upper bd prop}
Let $I$ be a height two Cohen-Macaulay ideal with
degree matrix
$$
A = \begin{bmatrix}
a_1 & b_1 & && *\\
d & a_2 & b_2 \\
&& \ddots & \ddots &  \\
* &  & & a_t & b_t
\end{bmatrix}.
$$
as in Remark \ref{rem-deg-matrix}
(but note the new variable $d$ in the $(1,2)$ spot if $t \geq 2$).
Then either one of the following conditions is sufficient to conclude that
\[
e(R/I) \leq \frac{1}{2} M_1M_2 - (M_1-m_1) - (M_2 - m_2).
\]
\begin{itemize}
\item[(i)] all of the entries of $A$ are $\geq 2$;
\item[(ii)] $a_1-2d+1 \geq 0$, provided $t \geq 2$.
\end{itemize}
\end{proposition}
\begin{proof}
We omit the details and leave the proof to the reader.
\end{proof}

\begin{example}
Consider the degree matrix
$$
B =
\left [
\begin{array}{ccccc}
2 & 2 & 2 & 2 \\
2 & 2 & 2 & 2 \\
1 & 1 & 1 & 1\end{array}
\right ]
$$
This comes, for example, from the ideal $I' = (x5, x^4y, x^2y3,y5)$.  It is
easy to check that $M_1' = m_1' = 5, M_2' = 7, m_2' = 6$, and
$e(R/I') = 17$.  If
we consider the
$2 \times 3$ submatrix $A$, as in the proof, we have $m_1 = M_1 = 4, m_2
= M_2 = 6$.  Then
$$
e(R/I')>\frac{1}{2} M_1'M_2' - (M_1'-m_1') - (M_2'-m_2') = 16.5.
$$
So we see that the conclusion of Proposition \ref{upper bd prop} does not
hold here, and indeed \linebreak $a_1-2d+1 = -1$.

In fact, it is not too difficult to show that a $t \times (t+1)$ degree
matrix consisting of 2's in the first $(t-1)$ rows and 1's in
the last row satisfies $m_1 = M_1 = 2t-1$, $m_2 = 2t$, $M_2 = 2t+1$,
and $e(R/I')
= 2t2-1$.  However, one checks that
\[
\frac{1}{2} M_2M_2 - (M_1-m_1) - (M_2-m_2) = 2t2 - \frac{3}{2}
\]
so this gives an example of any size that violates the bound of Proposition
\ref{upper bd prop}.
\end{example}

As mentioned in the introduction, it was shown by Huneke and Miller \cite{HM} that
if $R/I$ is
Cohen-Macaulay of codimension $c$ with a pure resolution then
\[
e(R/I) = \left ( \prod_{i=1}^c d_i \right )/ c! ,
\]
where of course $d_i = m_i = M_i$ is the shift in the $i$-th free module of the
pure resolution.  Corollary \ref{cor-intro-2} in the introduction may be seen as a
 converse to this result.

\begin{proof}[Proof of Corollary \ref{cor-intro-2}]
The claim follows from Theorem \ref{thm-intro-2} because $m_2 > m_1$ and $M_2 > M_1$.
\end{proof}

We end this section with a remark about the non Cohen-Macaulay case in
codimension two.

\begin{remark} \label{rem-naCM}

It is known that the lower bound of Conjecture \ref{conj} is
false in the non
Cohen-Macaulay case.  Even the weaker statement
\[
e(R/I) \geq \frac{1}{c!} \prod_{i=1}^c m_i
\]
(where $c$ is the codimension of $I$) is false.  An easy example is the
case  of two skew lines in
${\mathbb P}^3$.  In codimension two,  however, the analogous upper
bound is true
(\cite{roemer}).  Is there a different lower bound that is true?  One
natural guess is that
one might be able to replace $c!$ by some suitable integer $k$.  That
is, perhaps
\[
e(R/I) \geq \frac{1}{k} \prod_{i=1}^c m_i
\]
for suitable $k$.
    One immediately sees
that at the very least, we should assume that our ideals are unmixed.
For instance,
starting with any curve, adding  points does not change the
multiplicity but makes the
$m_i$ arbitrarily large.

The next observation is that we must assume that $I$ is reduced in
order to hope for a
lower bound of the type $\frac{1}{k} \prod_{i=1}^c m_i$.  Indeed,
consider ideals in
$k[x_0,x_1,x_2,x_3]$ of the form
\[
I = (x_0,x_1)^t + (F)
\]
where $F$ is a polynomial that is smooth along the line defined by
$(x_0,x_1)$, and $\deg
F \geq t+1$.  Then
$I$ defines an unmixed curve of multiplicity $t$ (cf.\ \cite{MPP}), and
one quickly sees that
$m_1 = t$ and $m_2 = t+1$.  Hence we would need $k \geq t+1$, which
can be made arbitrarily
large by choosing large $t$.

However, if we do assume that $I$ is reduced, there may be such a
bound.  Indeed,
experiments with {\tt Macaulay} (\cite{macaulay}) have not yet
yielded a counterexample to
the guess
\[
e(R/I) \geq \frac{1}{5} m_1 m_2
\]
at least among unmixed height two {\em reduced} non arithmetically
Cohen-Macaulay curves in
${\mathbb P}^3$.

\end{remark}


\section{Codimension three Gorenstein algebras}
\label{sec-3}

We now turn to height three Gorenstein ideals.  Our approach here is similar to
that of the alternative proofs given in the previous section. In
\cite{HS} the upper bound stated in Conjecture \ref{conj} was proved for such
ideals.  The lower bound was proved only when the number of generators is five
(or of course three).  In this section we will prove an improved version of
Conjecture
\ref{conj}, and as a consequence we will again (as in the codimension two case)
immediately obtain as a corollary that sharpness occurs (necessarily for both
bounds) if and only if the resolution is pure.

Let $I \subset R$ be a height three Gorenstein ideal.
The possible graded Betti numbers of such ideals were
described in \cite{BE} and in \cite{D}, and it was shown in \cite{GM5} that any
such set of graded Betti numbers in fact occurs for some reduced arithmetically
Gorenstein set of points in $\mathbb P3$.  In fact more was shown in
\cite{GM5}.

\begin{remark} \label{rem-Gor-const}
Suppose that $I$ has a minimal graded free resolution
\begin{equation} \label{desired resol}
0 \rightarrow R(-m_3) \rightarrow \bigoplus_{i=1}^{2t+1} R(-\be_i)
\rightarrow
\bigoplus_{i=1}^{2t+1} R(-\al_i) \rightarrow I \rightarrow 0
\end{equation}
where $\al_1 \leq \dots \leq \al_{2t+1}$, $\be_1 \leq \dots \leq \be_{2t+1}$ (this
is slightly different from the notation of \cite{GM5}), and $m_3 = M_3$.
It was shown in \cite{GM5} that

\begin{itemize}
\item \label{fact1} there exists a Cohen-Macaulay ideal
$J \subset R$ with minimal graded free resolution
\begin{equation} \label{resol of C}
0 \rightarrow \bigoplus_{i=1}^t R(-\be_i) \rightarrow \bigoplus_{i=1}^{t+1}
R(-\al_i) \rightarrow J \rightarrow 0.
\end{equation}
\item \label{fact2} there are homogeneous polynomials $f, g \in J$ of degree
 $m_1(I) =
\al_1$ and
$M_2(I) = \be_{2t+1} = m_3-\al_1$, respectively, such that $\tilde{J} := (f, g) : J$ does
not have any components in common with $J$ (i.e.\  $J$ and $\tilde{J}$ are geometrically
linked by $(f,
g)$), thus $(f, g) = J \cap \tilde{J}$.
\item \label{fact4} the ideal $I := J + \tilde{J}$ has the desired minimal free resolution
(\ref{desired resol}).  (Note that not all  Gorenstein ideals arise in this way; this only
says that
numerically for any set of graded Betti numbers this construction produces a Gorenstein
ideal with the given Betti numbers, but this is enough
for our purposes.)
\end{itemize}
\end{remark}

A famous result of Buchsbaum and Eisenbud in \cite{BE} says that without loss
of generality we may choose bases so that the middle map of the resolution
(\ref{desired resol}) is represented by a skew symmetric matrix, $M$,
and that the minimal 
generators of $I$ are then given by the maximal Pfaffians of that matrix.
However, we may represent the degree matrix $B$ corresponding to $M$,
much as we did
in (\ref{deg mat}), so that the entries are increasing from bottom to top and
from left to right.  Then the resulting degree matrix $B$ is {\em symmetric} about the {\em
non-main} diagonal. In particular, we have
\begin{equation} \label{eq-3-deg-mat}
B  =  \begin{bmatrix}
b_t & & & & & & *\\
a_t & \ddots \\
& \ddots & b_1 \\
& & a_1 & d \\
& & & a_1 & b_1 & & \\
& & && \ddots & \ddots &  \\
* & & & & & a_t & b_t \\
\end{bmatrix}
\end{equation}
where
\begin{equation}
A  :=  \begin{bmatrix}
a_1 & b_1 & &  *\\
& \ddots & \ddots &  \\
* & & a_t & b_t \\
\end{bmatrix}
\end{equation}
is the degree matrix of the ideal $J$ (that has been used to produce $I$).

Furthermore, comparing the resolutions (\ref{desired resol}) and (\ref{resol of C}) we obtain in
conjunction with the formulas (\ref{formulas}) that
\begin{equation} \label{3-formulas}
\begin{array}{rcccl}
m_1 & = & m_1 (J) & = & a_1 + \ldots + a_t \\
m_2 & = &  m_2 (J)& = & a_1 + \ldots + a_t + b_t = m_1 + b_t \\
m_3 & = & & & d + 2 (b_1 + \ldots + b_t).
\end{array}
\end{equation}

We now are ready to show our improvement of Conjecture \ref{conj} for Gorenstein ideals of
codimension three.

\begin{proof}[Proof of Theorem \ref{ht 3 gor bd}]
Our proof will be by induction on
the size of the degree matrix, $A$, of the Buchsbaum-Eisenbud matrix of $I$.
First, let $t = 1$, i.e.\ $I$ is a complete intersection. Let $m_1, y, M_1$ be the degrees of the
minimal generators of $I$. Then $m_1 \leq y \leq M_1$ and we get

\begin{eqnarray*}
\lefteqn{ m_1 m_2 m_3 + (M_3-M_2)2(M_2-m_2+ M_1-m_1) } \\[1ex]
& = & m_1 (m_1 + y) (m_1 + y + M_1) +  2 m_12 (M_1 - m_1) \\[1ex]
& = & m_1 (2 y + m_1 - y) (3 M_1 +m_1+ y - 2 M_1) +  2 m_12 (M_1 - m_1) \\[1ex]
& = & 6 m_1 y M_1 + m_1 \left [ (m_1 + y) (m_1 + y - 2M_1) + (m_1 - y) 3 M_1 + 2 m_1 (M_1 - m_1)      \right ] \\[1ex]
& = & 6 e(R/I) + m_1 \left [ (m_1 + y) (y - M_1) + 3 (m_1 - y)  M_1 + (m_1 -y)  (M_1 - m_1)      \right ] \\[1ex]
& \leq & 6 e(R/I)
\end{eqnarray*}
proving the lower bound.

The upper bound is shown similarly. We have
\begin{eqnarray*}
\lefteqn{M_1 M_2 M_3 - \frac{1}{2} M_3(M_2-m_2+M_1-m_1) } \\[1ex]
& = & M_1 (y + M_1) (m_1 + y + M_1) - (m_1 + y + M_1) (M_1 - m_1) \\[1ex]
& = & M_1 (2 y + M_1 - y) (3 m_1 + y + M_1 - 2m_1) - (m_1 + y + M_1)
(M_1 - m_1) \\[1ex]
& = & 6 m_1 y M_1 + M_1 \left [ (M_1 - y) (m_1 + y + M_1) + 2 y (y +
  M_1 - 2m_1) \right ] \\
& & - (m_1 + y + M_1) (M_1 - m_1) \\[1ex]
& = & 6 e(R/I) + 2 y M_1(y - m_1) + (m_1 + y + M_1) M_1
(M_1 - y) \\[1ex]
& &+  \left [2 y M_1 - (m_1 + y + M_1) \right ] (M_1 - m_1) \\[1ex]
& \geq & 6 e(R/I)
\end{eqnarray*}
because if $M_1 = 1$ then we must have $m_1 = y = M_1$ and the last
estimate becomes an equality. But if $M_1 \geq 2$ then we get
$2 y M_1 - (m_1 + y + M_1) \geq 2 y M_1 - (2y + M_1) \geq 0$
because for any two integers $k, l \geq 2$ we have $k l \geq k+l$.
The upper bound follows.

Now assume $t \geq 1$ and let $I'$ be the Gorenstein ideals whose degree matrix is
\begin{equation} \label{eq-3-ind-deg-mat}
B'  =  \begin{bmatrix}
b_{t+1} \\
a_{t+1} & b_t & & & & & & *\\
& a_t & \ddots \\
& & \ddots & b_1 \\
& & & a_1 & d \\
& & & & a_1 & b_1 & & \\
& & & && \ddots & \ddots &  \\
&  & & & & & a_t & b_t \\
& * & & & & &  & a_{t+1} & b_{t+1}
\end{bmatrix}
\end{equation}
and that has been produced using the Cohen-Macaulay ideal $J'$ (cf.\ Remark \ref{rem-Gor-const})
with degree matrix
$$
A' = \begin{bmatrix}
a_1 & b_1 & & & *\\
& \ddots & \ddots &  \\
 & & a_t & b_t \\
*   & & & a_{t+1} & b_{t+1}
\end{bmatrix}.
$$
To simplify notation, we set as in the codimension two case
$$
a := a_{t+1}, \quad b := b_{t+1}, \quad \mbox{and} \; c := b_t.
$$
Let $I$ be the Gorenstein ideal whose Buchsbaum-Eisenbud matrix is
obtained from the Buchs\-baum-Eisenbud matrix of $I'$ by stripping the
top and bottom rows, and the rightmost and leftmost columns such that
the degree matrix is $B$ as in (\ref{eq-3-deg-mat}).
Let $m_1, m_2, m_3,$ $M_1, M_2, M_3$ be the invariants of $I$ and let
$m'_1, m'_2, m'_3, M'_1, M'_2, M'_3$ be those of $I'$. Self-duality of
the resolution of $R/I$ provides
$$
M_1 = m_3 - m_2 \quad \mbox{and} \quad M_2 = m_3 - m_1.
$$
It follows that
\begin{equation}
M_2 - m_2 = M_1 - m_1 = m_3 - m_1 - m_2 = M_1 + M_2 - M_3.
\end{equation}
Thus, the formulas (\ref{3-formulas}) provide
\begin{equation} \label{eq-rel-I-I'}
\begin{array}{rclcrcl}
m'_1 & = & m_1 + a & \mbox{ } & M'_1 & = & M_1 + b + c -a \\
m'_2 & = & m_2 + a + b - c & \mbox{ } & M'_2 & = & M_2 + 2 b - a \\
m'_3 & = & m_3 + 2 b & \mbox{ } & M'_3 & = & M_3 + 2b.
\end{array}
\end{equation}
We also need the relation between the multiplicities of $R/I$ and $R/I'$.

\begin{claim} \label{rel bet G and G'}
$
\begin{array}{rcl}
e(R/I') & =  & e(R/I) + b (m_1 + a) (M_2 + b -a).
\end{array} $
\end{claim}

To see this, we may assume temporarily that $R/I$ and $R/I'$ have dimension one. Thus, the
ideals $J$ and $J'$ used to produced $I$ and $I'$ (as in Remark \ref{rem-Gor-const}) define
curves. Denote their arithmetic genera by $g$ and $g'$, respectively. As preparation, we first
relate the multiplicities of $R/I$ and $R/J$ and then the genera of $J$ and
$J'$.

Using the notation of Remark \ref{rem-Gor-const}, we have that $I = J + \tilde{J}$ and that $\fc
:= J \cap \tilde{J}$ is a complete intersection of type
$(m_1, M_2)$. Hence, we have graded isomorphisms (cf., e.g., \cite{N-gorliaison}, Lemma 3.5)
$$
K_{R/J} (4- m_1 - M_2) \cong \tilde{J}/(J \cap \tilde{J}) \cong (J + \tilde{J})/J = I/J
$$
where $K_{R/J}$ denotes the canonical module of $R/J$.
It follows that
\begin{equation} \label{rel bet I and J}
e(R/I) = (m_1 +M_2-4) \cdot e(R/J) - (2g-2).
\end{equation}
Numerically, we may assume that $J'$ is a basic double link of $J$ (cf.\ (\ref{eq-J'-J})), i.e.\
there is a complete intersection $(f, g)$ of type $(m_1+a, b)$ such that $J' = f R + g J$.
Hence using the formula for
the genus of a complete intersection (see for instance \cite{book}, page 36), Proposition 4.1(b)
in \cite{N-gorliaison} provides
\begin{equation*}
g' = g + \frac{1}{2} b (m_1 + a) (m_1 + a + b -4) + b \cdot e(R/J).
\end{equation*}
Therefore, using formula (\ref{rel bet I and J}) for $I'$ as well as (\ref{eq-rel-I-I'}) we
obtain
\begin{eqnarray*}
e(R/I') & = & (m'_1 + M_2' - 4) \cdot e(R/J') - (2 g' - 2) \\
& = & (m_1 + M_2 + 2 b - 4) \left [ e(R/J) + b (m_1 + a)  \right ] \\
& & - \left [ 2 g +  b (m_1 + a) (m_1 + a + b -4) + 2 b \cdot e(R/J) - 2 \right ] \\
& = & (m_1 +M_2-4) \cdot e(R/J) - (2g-2) + b (m_1 + a) (M_2 + b -a)  \\
& = & e(R/I) + b (m_1 + a) (M_2 + b -a),
\end{eqnarray*}
as claimed.

Now we are ready for the induction step. We assume that the bounds hold for $I$, and we prove them
for $I'$.  We  will use the above numbered equations without comment.

We begin with the  lower bound. We have to show that
\begin{equation} \label{eq-low-b}
e(R/I')  \geq  \frac{1}{6} m'_1 m'_2 m'_3 + \frac{1}{3} (m'_1)2 (m'_3 - m'_1 - m'_2).
\end{equation}
Unfortunately, we need some rather lengthy computation. We have
\begin{eqnarray*}
\lefteqn{ m'_1 m'_2 m'_3 + 2 (m'_1)2 (m'_3 - m'_1 - m'_2) } \\[1ex]
& = &  (m_1 + a) (m_2 + a + b - c) (m_3 + 2b)  +  2 (m_1 + a)2 (m_3 -
m_1 - m_2 - (2a - b -c)) \\[1ex]
& = & m_1 m_2 m_3 + 2 m_12 (m_3 - m_1 - m_2) \\
& &  + \left [a (m_2 + a + b - c)(m_3 + 2b) + m_1 (m_2 + a + b -c) 2 b
  + m_1 (a+ b -c) m_3  \right ] \\
& & + 2 \left [ a (2 m_1 + a)(m_3 - 2 m_1 - c) - (m_1 + a)2 (2a - b -
  c) \right ] \\[1ex]
& \leq &  6 e(R/I)  \\
& &  + \left [a (m_1 + a + b)(m_3 + 2b) + m_1 (m_1 + a + b) 2 b + m_1
  (a+ b -c) m_3  \right ] \\
& & + 2 \left [ a (2 m_1 + a)(m_3 - m_1 - m_2) - (m_1 + a)2 (2a - b -
  c) \right ]
\end{eqnarray*}
by the induction hypothesis and because of $m_2 - c = m_1$ (by \ref{3-formulas}).

Using Claim \ref{rel bet G and G'} and essentially ordering for $m_3$ we obtain
\begin{eqnarray*}
\lefteqn{ m'_1 m'_2 m'_3 + 2 (m'_1)2 (m'_3 - m'_1 - m'_2) } \\[1ex]
& \leq &  6 e(R/I')  - 6 b (m_1 + a) (m_3 - m_1 + b -a)  \\
& &  + \left [a (m_1 + a + b)(m_3 + 2b) + m_1 (m_1 + a + b) 2 b + m_1
  (a+ b -c) m_3  \right ] \\
& & + 2 \left [ a (2 m_1 + a)(m_3 - m_1 - m_2) - (m_1 + a)2 (2a - b -
  c) \right ] \\[1ex]
& = & 6 e(R/I')  \\
& & + m_3 \left [ -6b (m_1 + a) + a (m_1 + a + b) + (a+b-c) m_1 \right ]
+ 2a (2m_1 + a) (m_3 - 2m_1 - c) \\
& & - 6b (m_1 + a) (-m_1 + b -a) + a (m_1 +a + b) 2b + m_1 (m_1 + a + b) 2b
- 2 (m_1 + a)2 (2a - b -c) \\[1ex]
& = & 6 e(R/I')  \\
& & + m_3 \left [ m_1 (2a -5b -c) + a (a-5b) \right ] + 2a (2m_1 + a) (m_3 - 2m_1 - c) \\
& & - 6b (m_1 + a) (-m_1 + b -a) + 2b (m_1 + a) (m_1 +a + b) - 2 (m_1 + a)2 (2a - b -c) \\[1ex]
& = & 6 e(R/I')  \\
& & + m_3 \left [ m_1 (2a -b -c) + a (a-b) - 4b (m_1 + a) \right ] + 2a (2m_1 + a) (m_3 - 2m_1 - c) \\
& & +  2b (m_1 + a) (4 m_1 + 4 a - 2b) - 2 (m_1 + a)2 (2a - b -c) \\[1ex]
& = & 6 e(R/I')  \\
& & + m_3 \left [ m_1 (2a -b -c) + a (a-b) \right ] + 2a (2m_1 + a) (m_3 - 2m_1 - c) \\
& & +  4b (m_1 + a) (- m_3 + 2 m_1 + 2 a - b) - 2 (m_1 + a)2 (2a - b -c) \\[1ex]
& = &  6 e(R/I') + a (a-b) m_3 + (2a-b-c) \left [ m_1 m_3 + 4b (m_1 + a) - 2 (m_1 + a)2 \right ] \\
& & + (m_3 - 2m_1 -c) \left [ 2a (2m_1 + a) - 4b (m_1 + a) \right ]
\end{eqnarray*}
where we used $-m_3 + 2m_1 + 2a - b = -m_3 + 2m_1 + c + 2a - b -c$.
Observing that $a \leq b,\; a \leq c$, and $m_3 \geq m_1 + m_2 = 2 m_1
+ c$ we get
\begin{eqnarray*}
\lefteqn{ m'_1 m'_2 m'_3 + 2 (m'_1)2 (m'_3 - m'_1 - m'_2) } \\[1ex]
& \leq & 6 e(R/I') + (2a-b-c) \left [ m_1 m_3 + 4b (m_1 + a) - 2 (m_1 + a)2 \right ] \\
& & + (m_3 - 2m_1 -c) \left [ 4 m_1 (a-b) + 2a (a -2b) \right ] \\[1ex]
& \leq & 6 e(R/I') + (2a-b-c) \left [ m_1 m_3 + 4b (m_1 + a) - 2 (m_1 + a)2 \right ] \\[1ex]
& \leq & 6 e(R/I')
\end{eqnarray*}
because
$$
m_1 m_3 + 4b (m_1 + a) - 2 (m_1 + a)2 = m_1 (m_3  - 2 m_1 + 4 (b-a)) + 2a (2b - a) \geq 0.
$$
This completes the proof of the lower bound.

Turning to the upper bound,
we have to show that:
\begin{equation}
e(R/I') \leq \frac{1}{6} M'_1 M'_2 M'_3 - \frac{1}{6}
M'_3 (M'_1 + M'_2 - M'_3).
\end{equation}
To start with, we have:
\begin{eqnarray*}
\lefteqn{ M'_1 M'_2 M'_3 -  M'_3 (M'_1 + M'_2 - M'_3) } \\[1ex]
& = & (M_1 + b+c-a) (M_2 + 2b - a) ( M_3+ 2b) \\
& & -  (M_3 + 2b) (M_1 + M_2 - M_3 + b + c -2a) \\[1ex]
& = & M_1 M_2 M_3 -  M_3 (M_1 + M_2 - M_3) \\
& & + (M_1 +b+c-a) (2b-a) (M_3 + 2b) + (M_1 + b+c-a) M_2  2b +
(b+c-a) M_2 M_3 \\
& & - M_3 (b+c-2a) - 2b (M_1 + M_2 - M_3 + b+c-2a) \\[1ex]
& \geq & 6 e(R/I) \\
& & + (M_2 + b-a) (2b-a) (M_3 + 2b) + (M_2 + b-a) M_2 2b + (b+c-a) M_2 M_3 \\
& & - M_3 (b+c-2a) - 2b (2 M_2 - M_3 + b-2a) \\[1ex]
\end{eqnarray*}
by the induction hypothesis and because of $M_2 = M_1 + c$.

Using Claim \ref{rel bet G and G'} and essentially ordering for $M_3$
we obtain
\begin{eqnarray*}
\lefteqn{ M'_1 M'_2 M'_3 -  M'_3 (M'_1 + M'_2 - M'_3) } \\[1ex]
& \geq & 6 e(R/I') - 6 b (M_3 - M_2 + a) (M_2 + b - a) \\
& & + (M_2 + b-a) (2b-a) (M_3 + 2b) + (M_2 + b-a) M_2 2b + (b+c-2a) M_2 M_3 \\
& & - M_3 (b+c-2a) - 2b (2 M_2 - M_3 + b-2a) \\[1ex]
& = & 6 e(R/I') \\
& & + M_3 \left [ (M_2 + b-a) (2b -a - 6b) + (b+c-a) M_2 \right ] \\
& & + (M_2 + b-a) \left [ (2b-a) 2b + 2b M_2 + 6b (M_2 -a) \right ] \\
& & - M_3 (b+c-2a) - 2b (2 M_2 - M_3 + b-2a) \\[1ex]
& = & 6 e(R/I') \\
& & + M_3 \left [ (M_2 + b-a) (- 4b) + (b+c-2a) M_2 - a (b-a) \right ] \\
& & + (M_2 + b-a) \left [ 8b M_2 + (2b-4a) 2b   \right ] \\
& & - M_3 (b+c-2a) - 2b (2 M_2 - M_3 + b-2a) \\[1ex]
& = & 6 e(R/I') \\
& & + M_3 \left [ (b+c-2a) M_2 - a (b-a) \right ] \\
& & + (M_2 + b-a) 4 b \left [ 2 M_2 - M_3 + b-2a    \right ] \\
& & - M_3 (b+c-2a) - 2b (2 M_2 - M_3 + b-2a) \\[1ex]
& = & 6 e(R/I') \\
& & + M_3 \left [ (b+c-2a) (M_2 - 1)  - a (b-a) \right ] \\
& & + (M_2 + b-a- \frac{1}{2}) 4 b \left [ 2 M_2 - M_3 + b-2a    \right ].  \\
\end{eqnarray*}
Observe that $b \geq a$, $c \geq a$, and $M_2 \geq m_2 \geq c + 1 \geq
a + 1$. It follows that
\begin{eqnarray*}
b+c - 2a & \geq &  b - a, \\
M_2 - 1 & \geq & a,
\end{eqnarray*}
thus
\begin{eqnarray*}
(b+c-2a) (M_2 -1) & \geq & a (b - a).
\end{eqnarray*}
Since we also have
$$
2 M_2 - M_3 + b-2a = (M_1 + M_2 - M_3) + (b+ c-2a) \geq 0
$$
we obtain
$$
M'_1 M'_2 M'_3 -  M'_3 (M'_1 + M'_2 - M'_3) \geq 6 e(R/I')
$$
and the upper bound follows. \end{proof}

%

\begin{remark} \label{rem-hema}
In \cite{srinivasan}, Srinivasan proved, compared to Conjecture \ref{conj},
stronger bounds for Gorenstein ideals of arbitrary codimension, but with
quasi-pure resolutions. A resolution is quasi-pure if $m_i \geq
M_{i-1}$ for all 
$i$.  In case of a codimension three Gorenstein ideal $I$ her bounds are
\begin{equation*}
\frac{1}{6} m_1 M_2 M_3 \leq e(R/I) \leq \frac{1}{6} M_1 m_2 m_3.
\end{equation*} 
Note that the lower bound is not even true for arbitrary complete
intersections. For example, a complete intersection of type $(2, 2,
5)$ gives a counterexample. On the other hand, the upper bound is true
for all complete intersections and we wonder if it is true for all
Gorenstein ideals of codimension three. 
\end{remark}

The method of proof of Theorem \ref{ht 3 gor bd}
also provides an explicit formula for the multiplicity of a Gorenstein ideal in
terms of the degrees of the entries of its Buchsbaum-Eisenbud matrix.

\begin{proposition} \label{prop-def-formula}
Let $I$ be a homogeneous Gorenstein ideal of codimension three with $2 t +1$
minimal generators. Order its Buchsbaum-Eisenbud matrix such that its degree
matrix is $B$ as in (\ref{eq-3-deg-mat}). Then we have
$$
e(R/I) = \sum_{j=1}^t b_j \cdot (a_1 + \ldots + a_j) \cdot (d +
\sum_{i=1}^{j-1} 
(2 b_i - a_i) + b_j - a_j).
$$
\end{proposition}

\begin{proof}
This follows easily from Claim \ref{rel bet G and G'}. Indeed, let $t
= 1$. Then 
$I$ is a complete intersection with minimal generators of degree $a_1, b_1, d +
b_1 - a_1$ and the claim follows.  Let $t \geq 2$. Then Claim \ref{rel
  bet G and 
G'} provides the assertion using the
 formulas (\ref{3-formulas}) and $M_2 = m_3 - m_1$.
\end{proof}


\end{document}